\documentclass[12pt, reqno]{amsart}
\usepackage{times}
\usepackage{amsmath,amsthm, amssymb, latexsym}
\usepackage{thmtools}
\usepackage{mathrsfs}
\usepackage{graphicx}
\newtheorem{theorem}{Theorem}[section]

\newtheorem{lemma}{Lemma}[section]
\newtheorem{example}{Example}[section]

\theoremstyle{remark}
\newtheorem{remark}{Remark}[section]

\begin{document}
\title{Algebraic Methods in Difference Sets and Bent Functions}
\author{Pradipkumar H. Keskar and Priyanka Kumari}
\maketitle
\centerline{ Department of Mathematics,}
\centerline{Birla Intitute of Technology and Science, Pilani (Pilani Campus),}
\centerline{ Pilani 333031, India.}
\centerline{emails: keskar@pilani.bits-pilani.ac.in;}
\centerline{priyanka.kumari@pilani.bits-pilani.ac.in}
\begin{abstract}
 We provide some applications of a polynomial criterion for difference sets. These include counting the  difference sets with specified parameters in terms of  Hilbert functions, in particular a count of  bent functions. We also consider the  question about the bentness of certain Boolean functions introduced by Carlet when the $\mathcal{C}$-condition introduced by him doesn't hold.

\noindent Key Words :Hilbert functions, $\mathcal{C}$-condition, flat, difference set, bent functions

\noindent MSC 2010 : 05B10, 6E30,13P25, 94C10, 11T71.
\end{abstract}
\today
\section{Introduction} In \cite{KeP}, a criterion was developed for subsets of a finite abelian group to be difference sets with specified parameters. The criterion was in terms of polynomial conditions on their characteristic vectors. This opens up the possibilities of algebro-geometric methods in studying difference sets. Recall that for a finite group $G$ of order $v$ and integers $k, \lambda$, a subset $D$ of $G$ is called a difference set of $G$ with parameters $(v, k, \lambda)$, or $(v, k, \lambda)$-difference set of $G$,  if $|D|=k$ and $|\{(g_1, g_2) \in D \times D : g_1g_2^{-1}=g\}|=\lambda$ for any non-identity element $g \in G$.

This paper discusses its applications in the study of difference sets and  bent functions, which are cryptographically significant.  Recall that for an even positive integer $t > 2$, a Boolean function of $t$ variables is a bent function if and only if its support is a  difference set in $\left(\mathbb{Z}/2\mathbb{Z}\right)^t$ with parameters $\left(2^t, 2^{(t-1)} \pm 2^{(t-2)/2},2^{(t-2)} \pm 2^{(t-2)/2}\right)$(where signs are chosen consistently).  

In Section 2, we show that the number of difference sets in an abelian group is given by the affine Hilbert function of a certain ideal in a polynomial ring. As a special case, the same holds for the number of bent functions of an even dimension. It may be pointed out that the count of bent functions is an important unresolved issue and even the known bounds for it are quite weak, see \cite{Tok} for the details. It is hoped that Computer Algebra software like Macaulay2, Singular, which facilitate computation of  Hilbert functions, will play a role in enhancing our knowledge in this direction. 

In Section 3, the criterion of \cite{KeP} is reformulated to characterize those exchanges of elements of a difference set which again lead to a difference set. The formulation is in terms of certain values of a polynomial function.  In subsequent sections, this criterion is applied to establish non-bentness of an infinite familiy  of certain functions introduced by Carlet in \cite{Car}, which we now introduce.

Let $\mathbb{F}_2 = \{0, 1\}$ be the field with two elements and let $m$ be a positive integer and $t=2m$. For any $x=(x_1, \dots, x_m), y = (y_1, \dots, y_m) \in \mathbb{F}_2^m $, let $x \cdot y = \sum_{i=1}^m x_iy_i \in \mathbb{F}_2$. Let $L$ be an $\mathbb{F}_2$-subspace of $\mathbb{F}_2^m$, $L^{\perp} = \{y \in \mathbb{F}_2^m : x \cdot y =0 \text{ for all }x \in L\}$ be the orthogonal complement of $L$ and let $1_{L^{\perp}} : \mathbb{F}_2^m \to \mathbb{F}_2$ be defined by $1_{L^{\perp}}(x) =1$ if $x \in L^{\perp}$ and $1_{L^{\perp}}(x) =0$ otherwise. For a permutation $\pi$ of $\mathbb{F}_2^m$, consider the function $f_{(\pi, L)} : \mathbb{F}_2^t=\mathbb{F}_2^m \times \mathbb{F}_2^m \to \mathbb{F}_2$ be defined by
$$f_{(\pi, L)}(x, y)= x \cdot \pi(y) + 1_{L^{\perp}}(x).$$
 Carlet found a class of bent functions called $\mathcal{C}$-class of bent functions. For this purpose, he introduced $\mathcal{C}$-condition on $(\pi, L)$ thus : $(\pi, L)$ satisfies the $\mathcal{C}$-condition if and only if for any $a \in \mathbb{F}_2^m$, $\pi^{-1}(a+L)$ is a flat (i.e. an affine subspace) in $\mathbb{F}_2^m$. He then showed that $\mathcal{C}$-condition is sufficient for bentness of $f_{(\pi, L)}$, that is, if $(\pi, L)$ satisfies $\mathcal{C}$-condition then $f_{(\pi, L)}$ is a bent function. The class of bent functions obtained in this manner is called the $\mathcal{C}$-class of bent functions. The $\mathcal{C}$-condition was further explored in \cite{MaG}. But it is not known that failure of $\mathcal{C}$-condition by $(\pi, L)$ implies non-bentness of $f_{(\pi, L)}$. Thus we need another machinary to prove non-bentness of $f_{(\pi, L)}$. As supports of bent functions are difference sets, the results of \cite{KeP} become relevant.

In this paper, we consider $\pi$ defined by
$$\pi(x_1, \dots, x_m) = \left((x_1 + P(x_2, \dots, x_m)), x_2, \dots, x_m\right)$$
 for all $(x_1, \dots, x_m) \in \mathbb{F}_2^m$, where $P(X_2, \dots, X_m) \in \mathbb{F}_2[X_2, \dots, X_m]$. Thus $\pi$ is induced by an important type of polynomial automorphisms of $\mathbb{F}_2^m$, called
elementary automorphisms, a generating set of the so called tame automorphism group, see \cite{vdE}. For several classes of $(\pi, L)$, we decide when the $\mathcal{C}$-condition is satisfied and in some examples where it is not satisfied, we conclude the non-bentness of $f_{(\pi, L)}$. 

\section{ Counting the Difference Sets}

 Let $G=\prod_{l=1}^{t}\left(\frac{\mathbb{Z}}{n_l\mathbb{Z}}\right)
	$ be an abelian group and let $v=|G|$. For any $i_l \in \frac{\mathbb{Z}}{n_l\mathbb{Z}}; 1\le l \le t$, let $i_l^* \in \{0,1,\dots,n_l-1\}$ be such that $i_l=i_l^*+n_l\mathbb{Z}$. For any subset $T$ of $G$, $\alpha = (\alpha_g :g \in G) \in \mathbb{C}^v$ is called the point representation or characteristic vector of $T$ if $\alpha_g=1$ for $g \in T$ and $\alpha_g=0$ otherwise. Let $X_1, \dots, X_t$ be independent variables over $\mathbb{C}$ and let $A_g, g \in G$ be $v$ independent variables over $\mathbb{C}[X_1, \dots, X_t]$. Letting $X=(X_1, \dots, X_t)$ and $A=(A_g : g \in G)$, we define  $\Psi = \Psi(X,A) \in \mathbb{C}[X, A]$ by  
\begin{align*}
\Psi  =& \left(\sum_{(i_1,\dots,i_t) \in G}A_{i_1\cdots i_t}X_1^{i_1^*}\cdots X_t^{i_t^*}\right)\left(\sum_{(i_1,\dots,i_t) \in G}A_{i_1\cdots i_t}X_1^{n_1-i_1^*}\cdots X_t^{n_t-i_t^*}\right) \\
&- \lambda \left(\sum_{(i_1, \dots, i_t) \in G}X_1^{i_1^*}\cdots X_t^{i_t^*}\right) - (k-\lambda).
\end{align*}
Further  let  $U = \{\xi =(\xi_1, \dots, \xi_t) \in \mathbb{C}^t : \xi_l^{n_l}=1 \text{ for all }1 \le l \le t\}$ and $P_g(A) =A_g^2-A_g\in \mathbb{C}[A]$ for all $g \in G$.
In Theorem 3.2 of \cite{KeP}, we have given a polynomial criterion for $(v, k, \lambda)$ difference set in $G$ as follows :

{\bf Polynomial Criterion :} For $\alpha = \left(\alpha_g : g \in G\right) \in \mathbb{C}^v$,   $\alpha$ is a point representation of a $(v, k, \lambda)$ difference set in $G$ if and only if $\alpha$ satisfies the equations $P_g(A) =0\text{ for all } g \in G$, and $\Psi(\xi, A)=0$ for all $\xi = (\xi_1, \dots, \xi_t) \in U$.

As a consequence, $(v, k, \lambda)$ difference sets in $G$ are in one-to one correspondence with the points of the zero-dimensional affine algebraic set $V\left(\Psi(\xi, A), P_g(A):  \xi \in U, g \in G \right)=\{ \alpha \in \mathbb{C}^v : \Psi(\xi, \alpha)=0 \text{ for all } \xi \in U, P_g(\alpha)=0 \text{ for all }g \in G\}$. This brings us to the concept of affine Hilbert function of an ideal in a polynomial ring. To define this, let $R=k[Z_1, \dots, Z_r]$ be the polynomial ring in $r$ variables over a field $k$ and for a non-negative integer $s$,  let $R_{\le s}=\{f(Z_1, \dots, Z_r) \in R : \deg{f} \le s\}$. For an ideal $I$ of $R$, let $I_{\le s}=I \cap R_{\le s}$. Note that $R_{\le s}$ is a finite dimensional $k$-vector space and  $I_{\le s}$ is its subspace. We define the affine Hilbert function of $I$ as the integer valued function  ${}^aHF_I$ of non-negative integers such that  ${}^aHF_I(s)=\text{dim}_k\left(\frac{R_{\le s}}{I_{\le s}}\right)$. It turns out that there is a polynomial, called affine Hibert Polynomial of $I$ and denoted by ${}^aHP_I$, whose values coincide with the values of the affine Hilbert function of $I$ for large integer values of $s$. Letting $I(V)=\{ f(Z_1, \dots, Z_r) \in R : f(z_1, \dots, z_r)=0 \text{ for all } (z_1, \dots, z_r) \in V\}$ for $V \subset k^r$,  Exercise 11 on p. 475 of  \cite{CLO} shows that if $|V|$ is finite then  ${}^aHP_{I(V)}$ is the constant polynomial $|V|$. We now apply this to $V=V\left(\Psi(\xi, A), P_g(A) :  \xi \in U, g \in G \right)$ to get the following
\begin{theorem}
 The number of $(v, k, \lambda)$-difference sets in the group  $G=\prod_{l=1}^{t}\left(\frac{\mathbb{Z}}{n_l\mathbb{Z}}\right)$ is given by  ${}^aHP_J(s)$ for any $s \in \mathbb{C}$ where $J$ is the ideal of $\mathbb{C}[A]$ generated by $\{\Psi(\xi, A), P_g(A) : \xi \in U, g \in G$\}.
\end{theorem}

{\bf Proof: } In the light of the above discussion,  the proof will be complete if we show that $I(V(J)) = J$. By Strong  Hilbert Nullstellensatz (Theorem 6 on p. 176 of \cite{CLO}), $I(V(J))= \sqrt{J}$.  This reduces our work to showing that $\sqrt{J} =J$. This follows from Theorem 8.14 on p. 343 of \cite{BeW}, since $P_g(A)$ is square-free and $\mathbb{C}$ is perfect.

\begin{remark} In stead of the affine Hilbert Polynomial of $J$, we could use the Hilbert Polynomial of the homogenization $J^h$ of $J$ as well. In some computational software, Hilbert Polynomial of a homogeneous ideal in a  graded ring  is easier to deal with, hence we also give an alternate  formulation of the above theorem. Let $B$ be an indeterminate over $\mathbb{C}[A]$ and let $\mathbb{C}[A]^h = \mathbb{C}[A][B]$.
The homogenization  $J^h$ is the ideal of $\mathbb{C}[A]^h$ generated by $\{f^h : f \in J\}$, where for any $f \in \mathbb{C}[A]$, $f^h \in\mathbb{C}[A]^h$ is defined by $f^h(A, B)= B^{\text{deg}(f)}f\left(\frac{A_g}{B} : g \in G\right)$. To obtain a finite generating set of $J^h$, note that by Theorem 4 on p. 388 of \cite{CLO}, if $S$ is a Gro\"bner basis of $J$ then $\{f^h : f \in S\}$ is the Gro\"bner basis of $J^h$.  To define Hilbert Function $HF_{J^h}$ of $J^h$, for any non-negative integer $s$, consider the $k$-vector spaces 
$$\mathbb{C}[A]^h_s=\{f \in R^h : f=0 \text{ or $f$ is homogeneous of degree }s\}$$
and $J^h_s=J^h \cap \mathbb{C}[A]^h_s$. We define $HF_{J^h}(s)=\text{dim}_k\left(\frac{\mathbb{C}[A]^h_s}{J^h_s}\right)$. By Theorem 12 on p. 464 of \cite{CLO}, ${}^aHF_J(s)=HF_{J^h}(s)$ for all non-negative integers $s$. This allows us to replace  the affine Hilbert Function ${}^aHF_J$ by the Hilbert Function $HF_{J^h}(s)$.
\end{remark}
\begin{remark}
 The above theorem also gives a count of all bent functions of $t$ variables. Since the set of all  bent functions with supports of size $2^{(t-1)} + 2^{(t-2)/2}$ and the set of those with supports of size $2^{(t-1)} - 2^{(t-2)/2}$ are disjoint of same cardinality, the count of the bent functions in $t$ variables for an even $t$ is given by $2HF_{J^h}(s)$ where $(v, k, \lambda)=\left(2^t, 2^{(t-1)} + 2^{(t-2)/2},2^{(t-2)} + 2^{(t-2)/2}\right)$.
\end{remark}

 \section{A Difference Set Criterion}

 The polynomial criterion  developed in \cite{KeP} for $(v, k, \lambda)$ abelian difference sets imposes some restrictions on exchanges of elements of a $(v, k_1, \lambda_1)$ difference set to get another $(v, k_2, \lambda_2)$ difference set.  By introducing a complex valued polynomial $\Delta(D_1, D_2)(X_1, \dots, X_t)$,  we make these restrictions explicit, in terms of its certain values. Alternately, this criterion can also be phrased in the language of group characters, following Theorem 11.18 on p. 224 of \cite{MoP}.

 Let $G=\prod_{l=1}^{t}\left(\frac{\mathbb{Z}}{n_l\mathbb{Z}}\right)
	$ be an abelian group. For any $i_l \in \frac{\mathbb{Z}}{n_l\mathbb{Z}}; 0\le l \le t$, let $i_l^* \in \{0,1,\dots,n_l-1\}$ be such that $i_l=i_l^*+n_l\mathbb{Z}$.\\
	For any $T\subset G$, let $$\rho_G(T)(X_1,\dots,X_t)=\sum_{(i_1,\dots,i_t)\in T}X_1^{i_1^*}\cdots X_t^{i_t^*} \in \mathbb{C}[X_1,\dots,X_t].$$
	Let $U=\{(\xi_1,\dots,\xi_t)\in \mathbb{C}^t: \xi_l^{n_l}=1 \text{ for all } 1\le l\le t \}$. For any $(\xi_1, \dots, \xi_t) \in U$  and $(i_1, \dots, i_t) \in G$, we define $\xi_1^{i_1} \cdots \xi_t^{i_t} = \xi_1^{i_1^*}\cdots \xi_t^{i_t^*}$.

  Now let $v=|G|$ and $k, \lambda$ be non negative integers.
	For any $D\subset G$, let $D^{(-1)}=\{-d:d\in G \}$. In $(3.2*)$ of \cite{KeP} it was proved that
\begin{align*}
\nonumber & D\text{ is a }(v, k,\lambda)\text{ difference set in }G\text{ if and only if}\\
	&\rho_{G}(D)(\xi_1,\dots,\xi_t)\rho_{G}(D^{(-1)})(\xi_1,\dots,\xi_t)-\lambda \rho_{G}(G)(\xi_1,\dots,\xi_t)-(k-\lambda)=0\\
	& \text{for all } (\xi_1,\dots,\xi_t) \in U.
	\end{align*}

Note that for any $\xi \in \mathbb{C}$ with $|\xi|=1$, we have $\xi^{-1} = \bar{\xi}$, the complex conjugate of $\xi$. It follows that $\rho_G\left(D^{(-1)}\right)(\xi_1, \dots, \xi_t)$ is the complex conjugate of $\rho_G(D)(\xi_1, \dots, \xi_t)$ and hence we get
\begin{eqnarray}
\nonumber & D\text{ is a }(v, k,\lambda)\text{ difference set in }G\text{ if and only if}\\
\nonumber &|\rho_{G}(D)(\xi_1,\dots,\xi_t)|^2-\lambda \rho_{G}(G)(\xi_1,\dots,\xi_t)-(k-\lambda)=0\\
	& \text{for all } (\xi_1,\dots,\xi_t) \in U.
	\end{eqnarray}

	Now let us assume $n_i=2$ for all $i=1,\dots,t$ then $\rho_{G}(D)(\xi_1,\dots,\xi_t) \in \mathbb{R}$.
	Further $\rho_G(G)(\xi_1,\dots,\xi_t)=0$ if $\xi_i\ne 1$ for some $i\in \{1,\dots,t\}$,
	while $\rho_{G}(D)(\xi_1,\dots,\xi_t)=k$ and $\rho_G(G)(\xi_1,\dots,\xi_t)=v$ if $\xi_i= 1$ for all $i\in \{1,\dots,t\}$. This has the following consequence :
	
\begin{align}
\nonumber & \text{If }n_i=2\text{ for all }i=1,\dots,t\text{ then }D\text{ is a }(v, k, \lambda)\text{ difference set in }G\\
\nonumber & \text{ if and only if for any }(\xi_1,\dots,\xi_t) \in U\\
	& \left(\rho_G(D)(\xi_1,\dots,\xi_t)\right)= \begin{cases}
	k & \text{if } \xi_i=1 \text{ for all } i \in \{1,\dots,t\}\\	
	\pm \sqrt{k-\lambda} & \text{otherwise}.
	\end{cases}
\end{align}

Now suppose $D_1$ is a $(v, k_1, \lambda_1)$ difference set in $G$ and $D_2 \subset G$. Let
\begin{align*}
\Delta(D_1,D_2)(X_1, \dots, X_t) = & \rho_G(D_1 \setminus D_2)(X_1, \dots, X_t) \\
&- \rho_G(D_2 \setminus D_1)(X_1, \dots, X_t).
\end{align*}
Then we have
$$\Delta(D_1,D_2)(X_1, \dots, X_t) = \rho_G(D_1)(X_1, \dots, X_t) - \rho_G(D_2)(X_1, \dots, X_t)$$
and hence :
\begin{align}
\nonumber & \text{If }n_i=2\text{ for all }i=1,\dots,t\text{ and } D_1 \text{ is a} (v, k_1, \lambda_1) \text{ difference set in } G\\
\nonumber & \text{ then }D_2 \text{ is a }(v, k_2, \lambda_2)\text{ difference set in }G
 \text{ if and only if }\\
 \nonumber &\text{for any }\xi=(\xi_1,\dots,\xi_t) \in U \\
	&\Delta(D_1, D_2)(\xi) \in \begin{cases}
	\{k_1-k_2\} \\
 \text{if } \xi_i=1 \text{ for all } i \in \{1,\dots,t\};\\	
 & {}\\
	\{\sqrt{k_1-\lambda_1}-\sqrt{k_2-\lambda_2}, \sqrt{k_1-\lambda_1}+\sqrt{k_2-\lambda_2}\} \\
 \text{if }\rho_G(D_1)(\xi)=\sqrt{k_1-\lambda_1};\\
 & {}\\
 \{-\sqrt{k_1-\lambda_1}-\sqrt{k_2-\lambda_2}, -\sqrt{k_1-\lambda_1}+\sqrt{k_2-\lambda_2}\} \\
 \text{if }\rho_G(D_1)(\xi)=-\sqrt{k_1-\lambda_1}.
	\end{cases}
\end{align}

Moreover, if $(v, k_1, \lambda_1) = (v, k_2, \lambda_2)$ then
\begin{align}
\nonumber & D_2 \text{ is a }(v, k_1, \lambda_1)\text{ difference set in }G
 \text{ if and only if }\\
 \nonumber &\text{for any }\xi=(\xi_1,\dots,\xi_t) \in U
 \end{align}
 \begin{align}
	&\Delta(D_1, D_2)(\xi) \in \begin{cases}
	\{0\} \\
 \text{if } \xi_i=1 \text{ for all } i \in \{1,\dots,t\};\\	
 & {} \\
	\{0, 2\sqrt{k_1-\lambda_1}\} \\
 \text{if }\rho_G(D_1)(\xi)=\sqrt{k_1-\lambda_1};\\
 & {} \\
 \{-2\sqrt{k_1-\lambda_1}, 0\} \\
 \text{if }\rho_G(D_1)(\xi)=-\sqrt{k_1-\lambda_1}.
	\end{cases}
\end{align}

\section{The Analysis of $\mathcal{C}$-condition}

In the rest of the paper, we continue with the notation and terminology introduced in Section 1. Moreover we identify $\mathbb{F}_2^{r_1} \times \cdots \times \mathbb{F}_2^{r_u}$ with $\mathbb{F}_2^{r_1 + \cdots + r_u}$ in a natural way. Also for any integer $u \ge 0$, we denote by $\mathbf{0}_u$ the element of $\mathbb{F}_2^u$ whose all components are $0$ and by $\mathbf{1}_u$ the element of $\mathbb{F}_2^u$ whose all components are $1$. 

To search for examples when  $\mathcal{C}$-condition is not satisfied, we study the  $\mathcal{C}$-condition for  $(\pi, L)$.

 Since for any $x \in \mathbb{F}_2$, $x^2 =x$, by reducing all the exponents of variables mod 2, without loss of generality we can assume
$$P(X_2,\dots, X_m)=\sum_{ \ell =0}^{m-1} \ \sum_{2 \le i_1 < i_2 < \dots < i_{\ell} \le m} \alpha_{i_1\cdots i_\ell}X_{i_1}\cdots X_{i_\ell}$$
where $\alpha_{i_1\cdots i_\ell} \in \mathbb{F}_2$ for all $2 \le i_1 < i_2 < \dots < i_{\ell} \le m$.

To search for examples when  $\mathcal{C}$-condition is not sarisfied, the search space is provided by the following :

\begin{theorem} Let $m \ge 2$ and $s \in [1, m]$ be integers. \newline
(A) Let $L=\{(x_1,\dots,x_s,\mathbf{0}_{m-s}):x_i \in \mathbb{F}_2, 1\le i \le s \}$ be a linear subspace of $\mathbb{F}_2^m$.  Then $\mathcal{C}$-condition is satisfied by $(\pi, L)$. \newline
(B) Let $L=\{(\mathbf{0}_s,x_{s+1},\dots,x_{m}): x_i \in \mathbb{F}_2, s+1\le i \le m \}$ be a linear subspace of $\mathbb{F}_2^m$. If $\alpha_{i_1\cdots i_\ell}=0$ for {\it all} $(i_1,\dots, i_\ell)$ such that $|\{i_j:i_j>s\}|\geq 2$ then $\mathcal{C}$-condition is satisfied  by $(\pi, L)$. \newline
(C) Let $L=\{(\mathbf{0}_s,x_{s+1},\dots,x_{m}): x_i \in \mathbb{F}_2, s+1\le i \le m \}$ be linear subspace of $\mathbb{F}_2^m$.
If $\alpha_{i_1\cdots i_\ell}=1$ for some $(i_1,\dots, i_\ell)$ such that $|\{i_j:i_j>s\}|\ge 2$ then $\mathcal{C}$-condition is not satisfied by $(\pi, L)$.

Moreover in (A) and (B), $f_{(\pi, L)}$ is a $\mathcal{C}$-class bent function.
\end{theorem}

\noindent{\bf Proof :} The proof is based on the following observation :
\begin{flalign}
& \text{A  nonempty subset }F \subset \mathbb{F}_2^m \text{F is a flat} \notag\\
& \Leftrightarrow F - b \text{ is a subspace of }\mathbb{F}_2^m \text{ for some }b \in F \notag\\
& \Leftrightarrow F - b \text{ is a subspace of }\mathbb{F}_2^m \text{ for all }b \in F. &&
\end{flalign}

We will apply this when $F=\pi^{-1}(a+L)$ and $b =\pi^{-1}(a)$ where $a \in  \mathbb{F}_2^m$.

(A) For any $a=(a_1,\dots,a_m)\in \mathbb{F}_2^m$ we see that $a+L=a^*+L$ where $a_j^*=a_j$ if $j> s$ and $a_j^*=0$ otherwise. Thus we can assume, without loss of generality, that $a_j=0$ for $j\le s$. Therefore $a+L=\{(x_1,\dots,x_s,a_{s+1},\dots,a_m):x_1,\dots,x_s \in \mathbb{F}_2 \}$. It is enough to show that $\pi^{-1}(a+L)-\pi^{-1}(a)$ is a subspace of $\mathbb{F}_2^m$ for all $a=(\mathbf{0}_s,a_{s+1},\dots,a_m)\in \mathbb{F}_2^m$.\\
	Clearly, $\pi^{-1}(a+L) -\pi^{-1}(a) \subset L$ and they have same cardinality. 
 Since $|L|$ is finite, we get $\pi^{-1}(a+L)-\pi^{-1}(a) = L$. 	

Thus $\pi^{-1}(a+L)-\pi^{-1}(a)$ is a subspace for all  $a=(\mathbf{0}_s,a_{s+1},\dots,a_m)\in \mathbb{F}_2^m.$ Hence $(\pi,L)$ satisfies the $\mathcal{C}$-condition.\\

\noindent (B) For any $a=(a_1,\dots,a_m)\in \mathbb{F}_2^m$ we see that $a+L=a^*+L$ where $a_j^*=a_j$ if $j\le s$ and $a_j^*=0$ otherwise. Thus we can assume, without loss of generality, that $a_j=0$ for $j>s$. Therefore $a+L=\{(a_1,\dots,a_s,x_{s+1},\dots,x_m):x_{s+1},\dots,x_m\in \mathbb{F}_2 \}$. Consequently it is enough to show that $\pi^{-1}(a+L)-\pi^{-1}(a)$ is a subspace of $\mathbb{F}_2^m$ for all $a=(a_1,\dots,a_s,\mathbf{0}_{m-s})\in \mathbb{F}_2^m$.

For any $(i_1,\dots, i_\ell)$, if $\alpha_{i_1\cdots i_\ell}=1$ then $\prod_{i \in \{i_j : i_j \le s\}} a_i$ is a constant in $\mathbb{F}_2$, and $|\{i_j : i_j >s \}| \le 1$. So for any $a=(a_1,\dots,a_s,\mathbf{0}_{m-s})\in \mathbb{F}_2^m$,
{\footnotesize $$\pi^{-1}(a+L)-\pi^{-1}(a)=\{(l(x_{s+1},\dots,x_m),\mathbf{0}_{s-1},x_{s+1},\dots,x_m):x_{s+1},\dots,x_m\in \mathbb{F}_2 \}$$}
 where $l(X_{s+1},\dots,X_m) \in \mathbb{F}_2[X_{s+1},\dots,X_m] $ is a polynomial of degree $\le1$. \\
Now for any $u,v \in \pi^{-1}(a+L)-\pi^{-1}(a)$ and $\alpha,\beta \in \mathbb{F}_2$, where
\begin{align*}
&u=(l(x_{s+1},\dots,x_m),\mathbf{0}_{s-1},x_{s+1},\dots,x_m) \text{ and }\\ &v=(l(x_{s+1}^*,\dots,x_m^*),\mathbf{0}_{s-1},x_{s+1}^*,\dots,x_m^*),
\end{align*}
{\footnotesize
\begin{align*}
&\alpha u+\beta v\\
 &=\alpha(l(x_{s+1},\dots,x_m),\mathbf{0}_{s-1},x_{s+1},\dots,x_m) +\beta(l(x_{s+1}^*,\dots,x_m^*),\mathbf{0}_{s-1},x_{s+1}^*,\dots,x_m^*)\\
&=(\alpha l(x_{s+1},\dots,x_m)+\beta l(x_{s+1}^*,\dots,x_m^*),\mathbf{0}_{s-1},\alpha x_{s+1}+\beta x_{s+1}^*,\dots,\alpha x_m+\beta x_m^*)
\end{align*}}

Now $\pi^{-1}(a+\mathbf{0}_m)-\pi^{-1}(a)=\mathbf{0}_m$. Therefore $l(X_{s+1},\dots,X_m)$ is a polynomial with no constant term. Hence $l(X_{s+1},\dots,X_m)$ is a linear transformation. Then
{\footnotesize \begin{align*}
&\alpha u+\beta v\\
&= (l(\alpha x_{s+1}+\beta x_{s+1}^*,\dots,\alpha x_m+\beta x_m^*),\mathbf{0}_{s-1},\alpha x_{s+1}+\beta x_{s+1}^*,\dots,\alpha x_m+\beta x_m^*),
\end{align*}}
therefore $\alpha u+\beta v \in \pi^{-1}(a+L)-\pi^{-1}(a).$\\
As a consequence, $\pi^{-1}(a+L)-\pi^{-1}(a)$ is a subspace of $\mathbb{F}_2^m$ for all $a=(a_1,\dots,a_s,\mathbf{0}_{m-s})\in \mathbb{F}_2^m$ and hence $(\pi,L)$ satisfy $\mathcal{C}$-condition.\\

\noindent (C) We classify nonzero terms of $P(X_2, \dots, X_m)$ in two types.
\begin{itemize}
\item[Type 1 :] corresponding to $(i_1,\dots, i_\ell)$ such that $|\{i_j:i_j>s\}|\ge 2$,
\item[Type 2 :]  corresponding to $(i_1,\dots, i_\ell)$ such that $|\{i_j:i_j>s\}| < 2$.
\end{itemize}
Let $T_1=X_{i^{*}_1} \cdots X_{i^{*}_{\ell^{*}}}$ be minimal among all nonzero terms of $P(X_2, \dots, X_m)$ of Type 1 with the divisibility partial order. Hence for every nonzero term  $T \ne T_1$  of Type 1 of $P(X_2, \dots, X_m)$ corresponding to $(i_1, \dots, i_{\ell})$ , there exists $1 \le j \le \ell$ such that $i_j \not\in \{i^{*}_1, \dots, i^{*}_{\ell^{*}}\}$, and therefore $T$ is divisible by $X_{i_j}$. For any $a = (a_1, \dots, a_m) \in \mathbb{F}_2^m$ we see that $a+ L = a^{*} +L$ where $a^{*}_j = a_j$ if $j \le s$ and $a^{*}_j =0$ otherwise. Thus we can assume, without loss, that $a_j =0$ for $j >s$. Therefore $a+L = \{ (a_1, \dots, a_s, x_{s+1}, \dots, x_m) : x_{s+1}, \dots, x_m \in \mathbb{F}_2\}$. In view of $(5)$, we want to show that
$\pi^{-1}(a+L) - \pi^{-1}(a)$ is not a subspace of $\mathbb{F}_2^m$ for some $a = (a_1, \dots, a_s, \mathbf{0}_{m-s}) \in \mathbb{F}_2^m$.

Let $a_j = 1$ for all $j = i^{*}_u \le s$ and $a_j =0$ for any $j \in \{1, 2, \dots, s\} \setminus \{i^{*}_1, \dots, i^{*}_{\ell^{*}}\}$. Since any term of $P(X_2, \dots, X_m)$ of Type 1 except $T_1$ is divisible by $X_{i_j}$ for some $i_j \not\in \{i^{*}_1, \dots, i^{*}_{\ell^{*}}\}$, in addition if we let $x_j=0$ for all $j \in \{s+1, \dots, m\} \setminus \{i^{*}_1, \dots, i^{*}_{\ell^{*}}\}$ then
\begin{align*}
& \pi^{-1}(a_1, \dots, a_s, x_{s+1}, \dots, x_m)\\
&= \left( \left(a_1 + l(x_{s+1}, \dots, x_m) + \prod_{r=s_0}^{\ell^{*}}x_{i^{*}_r}\right), a_2, \dots, a_s, x_{s+1}, \dots, x_m\right)
\end{align*}
 where $s_0 = \min\{j \in \{1, \dots, \ell^{*}\} : i^{*}_j > s \}$ and $l(X_{s+1}, \dots, X_m)$ is a polynomial of degree $\le 1$ coming from terms of Type 2. Therefore, in this case,
\begin{align*}
&\pi^{-1}(a+L)-\pi^{-1}(a)\\
&=\resizebox{.95\hsize}{!}{$\Bigg\{ \left(\left( l(x_{s+1}, \dots, x_m)-l(\mathbf{0}_{m-s}) + \prod_{r=s_0}^{\ell^{*}} x_{i^{*}_r}\right), \mathbf{0}_{s-1}, x_{s+1}, \dots, x_m\right)
	: x_{s+1}, \dots, x_m \in \mathbb{F}_2\Bigg\}.$}
\end{align*}

 For $s_0 \le j \le \ell^{*}$, let $e_{i^{*}_j}=( x_{s+1}, \dots, x_m)$ be such that $x_{i^{*}_j}=1$ and $x_i =0$ for $i \ne i^{*}_j$ and let $f_{i^{*}_j}=\pi^{-1}(a+(\mathbf{0}_s, e_{i^{*}_j})) - \pi^{-1}(a)$. Then $f_{i^{*}_j}= \left(l(e_{i^{*}_j})-l(\mathbf{0}_{m-s}),\mathbf{0}_{s-1}, e_{i^{*}_j}\right)$, as $x_{i^{*}_u}=0$ for $u \ne j$. Now $f_{i^{*}_j} \in \pi^{-1}(a+L) -\pi^{-1}(a)$ for all $s_0 \le j \le \ell^{*}$.  On the other hand, $\sum_{j=s_0}^{\ell^{*}} f_{i^{*}_j} = \left(l(\sum_{j=s_0}^{\ell^{*}} e_{i^{*}_j})-l(\mathbf{0}_{m-s}),\mathbf{0}_{s-1}, \sum_{j=s_0}^{\ell^{*}} e_{i^{*}_j}\right)$, as $l(X_{s+1}, \dots, X_m) -l(\mathbf{0}_{m-s})$ is a homogeneous linear polynomial. Denoting $\sum_{j=s_0}^{\ell^{*}} e_{i^{*}_j}$ by $(y_{s+1}, \dots, y_m)$, we have $\sum_{j=s_0}^{\ell^{*}} f_{i^{*}_j} = \left(l(y_{s+1}, \dots, y_m)-l(\mathbf{0}_{m-s}),\mathbf{0}_{s-1}, y_{s+1}, \dots, y_m\right)$. Since $\prod_{r=s_0}^{\ell^{*}}y_{i^{*}_r} =1$, we see that $\sum_{j=s_0}^{\ell^{*}} f_{i^{*}_j} \not\in \pi^{-1}(a+L) - \pi^{-1}(a)$.

 As a consequence, $\pi^{-1}(a+L) - \pi^{-1}(a)$ is not a subspace of $\mathbb{F}_2^m$ and hence $(\pi, L)$ does not satisfy $\mathcal{C}$-condition. \qed


 \section{ Non-bentness of an Infinite Family}
The violation of $\mathcal{C}$-condition by $(\pi, L)$ is not sufficient to show $f_{(\pi, L)}$ is not bent. Using the adaptation of difference set criterion from Section 3, we will show the non-bentness of $f_{(\pi, L)}$ for several $(\pi, L)$ in every even dimension. We require the following

\begin{lemma} Let $m \ge 3$ and $1\le s \le m-2$ be integers. Then
$$\sum_{(x_{s+1},\dots,x_m) \in \mathbb{F}_2^{m-s}} (-1)^{\left(\sum_{i={s+1}}^{m}x_i+\prod_{i={s+1}}^{m}(x_i+1)\right)}=-2.$$
\end{lemma}

\noindent\textbf{Proof:} More generally, we will prove : for any $j = s, s+1, \cdots, m-1$,
\begin{equation*}
\sum_{(x_{j+1},\dots,x_m) \in \mathbb{F}_2^{m-j}} (-1)^{\left(\sum_{i={j+1}}^{m}x_i+\prod_{i={j+1}}^{m}(x_i+1)\right)}=-2. \tag{\dag}
\end{equation*}
We prove $(\dag)$ by induction on $u = m-j$. If $u=1$, we have $j = m-1$. Since
$$\sum_{x_m \in \mathbb{F}_2} (-1)^{(1+(x_m+1))}(-1)^{x_m} =2,$$
$(\dag)$  holds for $u=1$.

Assume $(\dag)$ for $u = \nu$ where $\nu \le m-2$ and let $u = \nu+1$, that is, $j = m-\nu-1$. Now
\begin{align*}
&\sum_{(x_{j+1},\dots,x_m) \in \mathbb{F}_2^{m-j}} (-1)^{\left(\sum_{i={j+1}}^{m}x_i+\prod_{i={j+1}}^{m}(x_i+1)\right)}\\
&=\sum_{(x_{m-\nu},\dots,x_m) \in \mathbb{F}_2^{\nu+1}} (-1)^{\left(\sum_{i={m-\nu}}^{m}x_i+\prod_{i={m-\nu}}^{m}(x_i+1)\right)}\\
&=-\sum_{(x_{m-\nu},\dots,x_m) \in \mathbb{F}_2^{\nu+1}} \Big((-1)^{(1+\prod_{i={m-(\nu-1)}}^{m}(x_i+1))}\Big)^{(x_{m-\nu}+1)}(-1)^{(\sum_{i={m-(\nu-1)}}^{m}x_i)}\\
&= -\sum_{(x_{(m-\nu)+1},\dots,x_m) \in \mathbb{F}_2^{\nu}} (-1)^{(1+\prod_{i={(m-\nu)+1}}^{m}(x_i+1))}(-1)^{(\sum_{i={(m-\nu)+1}}^{m}x_i)}\\
&-\sum_{(x_{(m-\nu)+1},\dots,x_m) \in \mathbb{F}_2^{\nu}}(-1)^{(\sum_{i={(m-\nu)+1}}^{m}x_i)}
\end{align*}
Since $\sum_{x_m \in \mathbb{F}_2} (-1)^{x_m}=0$,
$$\sum_{(x_{(m-\nu)+1},\dots,x_m) \in \mathbb{F}_2^{\nu}}(-1)^{(\sum_{i={(m-\nu)+1)}}^{m}x_i)} = \prod_{i=(m-\nu)+1}^m \Big(\sum_{x_i \in \mathbb{F}_2} (-1)^{x_i}\Big)=0.$$
Moreover as $(\dag)$ holds for $u = \nu$, that is $j = m-\nu$, we have
$$\sum_{(x_{(m-\nu)+1},\dots,x_m) \in \mathbb{F}_2^{\nu}} (-1)^{(\prod_{i={(m-\nu)}}^{m}(x_i+1))}(-1)^{(\sum_{i={(m-\nu)+1}}^{m}x_i)}=-2.$$
As a consequence, $(\dag)$ holds for $u = \nu+1$. This completes the proof. \qed

Now we come to the main result.
	\begin{theorem} Let $m \ge 3$ and $1\le r\le s \le m-2$ be integers. Further let $L=\{(\mathbf{0}_s,x_{s+1},\dots,x_m): x_i \in \mathbb{F}_2, s+1\le i \le m  \}$ be an $m-s$ dimensional linear subspace of $\mathbb{F}_2^m$ and $\pi(x)=\left((x_1+\prod_{i=r+1}^{m}x_i),x_2,\dots,x_m\right)$ be a permutation of $\mathbb{F}_2^m$. Then $f_{(\pi,L)}:\mathbb{F}_2^{2m} \to \mathbb{F}_2$
		is not a bent function.
\end{theorem}
\noindent\textbf{Proof:}
Since $L=\{(\mathbf{0}_s,x_{s+1},\dots,x_m): x_i \in \mathbb{F}_2, s+1\le i \le m  \}$\\
we have $ L^{\perp}= \{(x_1,\dots,x_s,\mathbf{0}_{m-s}): x_i \in \mathbb{F}_2, 1\le i\le s \}$.\\
Also for any $y=(y_1,\dots,y_m) \in \mathbb{F}_2^m$, $\pi^{-1}(y)=\left(\left(y_1+\prod_{i=r+1}^{m}y_i\right),y_2,\dots,y_m\right)$. Therefore
\begin{align*}
f_{(\pi,L)}(x,y)&= \sum_{i=1}^{m}x_iy_i+x_1\prod_{i=r+1}^{m}y_i+\prod_{i={s+1}}^{m}(x_i+1)\\
&= f(x,y)+x_1\prod_{i=r+1}^{m}y_i
\end{align*}
where $f(x,y)=\sum_{i=1}^{m}x_iy_i+\prod_{i={s+1}}^{m}(x_i+1)$ is a $\mathcal{M}$-class bent function in $2m$ variables, see p. 90 of \cite{Dil}.\\
Let $D_{(\pi,L)}$ and $D$ denote the supports of $f_{(\pi,L)}$ and $f$ respectively.
Then
$$D \setminus D_{(\pi,L)}= \{(x,y)\in D:x_1\prod_{i=r+1}^{m}y_i=1 \}.$$
We know
$$x_1\prod_{i=r+1}^{m}y_i=1 \iff  x_1=y_{r+1}=\cdots =y_m=1.$$
Therefore
\begin{align*}
&(x,y)\in D \text{ and } x_1\prod_{i=r+1}^{m}y_i=1\\
& \iff y_1+\sum_{i=2}^{r}x_iy_i+\sum_{i=r+1}^{m}x_i+\prod_{i={s+1}}^{m}(x_i+1)=1 \text{ and }\\
& x_1=y_{r+1}=\cdots =y_m=1\\
& \iff y_1=1 +\sum_{i=2}^{r}x_iy_i+ \sum_{i=r+1}^{m}x_i+\prod_{i={s+1}}^{m}(x_i+1) \text{ and }\\
&x_1=y_{r+1}=\cdots =y_m=1.
\end{align*}

\noindent Consequently
\begin{align*}
& D \setminus D_{(\pi,L)} \\
&=\{(1,x_2,\dots,x_m,1 +\sum_{i=2}^{r}x_iy_i+ \sum_{i=r+1}^{m}x_i+\prod_{i={s+1}}^{m}(x_i+1),y_2,\dots,y_r,\mathbf{1}_{m-r}):\\
&x_2,\dots,x_m,y_2,\dots,y_r \in \mathbb{F}_2\}
\end{align*}
 and hence
$$|D\setminus D_{(\pi,L)}|=2^{m+r-2}.$$

Now, if $\bar{D}$ denotes the complement of $D$ in $\mathbb{F}_2^{2m}$,
\begin{align*}
&D_{(\pi,L)}\setminus D\\
& = \{(x,y)\in \bar{D}: x_1\prod_{i=r+1}^{m}y_i=1 \}\\
&=\{(1,x_2,\dots,x_m,\sum_{i=2}^{r}x_iy_i+ \sum_{i=r+1}^{m}x_i+\prod_{i={s+1}}^{m}(x_i+1),y_2,\dots,y_r,\mathbf{1}_{m-r}):\\
&x_2,\dots,x_m,y_2,\dots,y_r \in \mathbb{F}_2 \}
\end{align*}
and hence $$|D_{(\pi,L)}\setminus D|=2^{m+r-2}.$$\\
Let $U=\{1,-1 \}^{2m}$. Then for any $(\xi,\eta) \in U$ we have
\begin{align*}
&\Delta(D, D_{(\pi,L)})(\xi, \eta)=\xi_1\eta_{r+1}\cdots \eta_m(\eta_1-1)\times\\
&\sum_{(x_2,\dots,x_m,y_2,\dots,y_r) \in \mathbb{F}_2^{m+r-2}}\xi_2^{x_2}\cdots \xi_m^{x_m}\eta_1^{\left(\sum_{i=2}^{r}x_iy_i+ \sum_{i=r+1}^{m}x_i+\prod_{i={s+1}}^{m}(x_i+1)\right)}\eta_2^{y_2}\cdots \eta_r^{y_r}.
\end{align*}
Henceforth let $\xi_1=\cdots =\xi_r=1, \xi_{r+1}=\cdots=\xi_s=-1, \xi_{s+1}=\cdots=\xi_m=1,\eta_1=-1, \eta_2=\cdots=\eta_m=1$. Further let
\begin{align*}
\Lambda_1 &= \sum_{(x_2,\dots,x_r,y_2,\dots,y_r) \in \mathbb{F}_2^{2r-2}}(-1)^{\left(\sum_{i=2}^{r}x_iy_i\right)} \text{ and }\\
\Lambda_2 &= \sum_{(x_{s+1},\dots,x_m) \in \mathbb{F}_2^{m-s}}(-1)^{\left(\sum_{i={s+1}}^{m}x_i+\prod_{i={s+1}}^{m}(x_i+1)\right)}.
\end{align*}
Then
\begin{align*}
\Delta(D, D_{(\pi,L)})(\xi, \eta) &= -2 \sum_{(x_{r+1}, \dots, x_s) \in \mathbb{F}_2^{s-r}} \left(\Lambda_1 \Lambda_2\right)\\
&= -2^{s-r+1}\left(\Lambda_1 \Lambda_2\right).
\end{align*}
Now
$\Lambda_1=\left(\sum_{(x_2,y_2)\in \mathbb{F}_2^2}(-1)^{x_2y_2}\right)\cdots \left(\sum_{(x_r,y_r)\in \mathbb{F}_2^2}(-1)^{x_ry_r}\right)$ and \newline$\sum_{(x_i,y_i)\in \mathbb{F}_2^2}(-1)^{x_iy_i}=2$ for any $i=2,\dots,r$. Hence by Lemma 5.1,
$$\Delta(D, D_{(\pi,L)})(\xi, \eta)=-2^{s}\Lambda_2 =2^{s+1}.$$

 Since $|D \setminus D_{(\pi,L)}| = |D_{(\pi,L)} \setminus D| = 2^{m-1}$, we have $|D|=|D_{(\pi,L)}|$. If $f_{(\pi,L)}$ is a bent function, then $D_{(\pi,L)}$ is a difference set. Thus by Ryser's condition (see Section 3 of \cite{KeP}), it follows that parameters $(v, k, \lambda)$ of $D_{(\pi,L)}$ are same as those of $D$, hence $k -\lambda = 2^{t-2}=2^{2(m-1)}$. Consequently, by $(4)$ of Section 3, for any $(\xi, \eta) \in U$,
$\Delta(D, D_{(\pi,L)}) \in \{0, \pm 2\sqrt{(k-\lambda)}\} = \{0, \pm 2^m\}$.

But we have found $(\xi, \eta) \in U$ such that $\Delta(D, D_{(\pi,L)})=2^{s+1} \not\in \{0, \pm 2^m\}$ for any $1\le s \le m-2$. Therefore  $f_{(\pi,L)}$ is not bent function.
\qed

\begin{remark}
The alternate tools for proving the non-bentness of $f_{(\pi, L)}$ could have been Theorem on p. 94 of \cite{Car} or Proposition 1 on p. 398 of \cite{Kol}. As far as \cite{Car} is concerned, checking the required conditions is tedious. On the other hand, to prove non-bentness of $f_{(\pi, L)}$ using  \cite{Kol},it is sufficient to have either (a) the Hamming distance $d(f, f_{(\pi, L)}) < 2^m$ or (b) $d(f, f_{(\pi, L)}) = 2^m$ and either support $\mathfrak{A}$ of $f+ f_{(\pi, L)}$ is not a flat or the restriction of $f$ to $\mathfrak{A}$ is not an affine function. Now $d\left({f,f_{(\pi,L)}}\right)=|D\setminus D_{(\pi,L)}|+|D_{(\pi,L)}\setminus D|=2^{m+r-1}$. So if $r>1$ then Proposition 1 of \cite{Kol} doesn't help, though it suffices for $r=1$.
\end{remark}

\section{ Computational Results}
In this section, we report some more $(\pi, L)$ such that $f_{(\pi, L)}$ is not bent. This was established through implementations of methods of previous sections as well as \cite{KeP} using programs in C++ language.
In Section 5, the polynomial $P(X_2, \cdots, X_m)$ contained only one term of Type 1 (as described in case (C) of Theorem 4.1). Examples 6.1 and 6.2 contain more than one term of Type 1.
\begin{example}
Let $L=\{(0,0,x_3,x_4): x_3,x_4\in \mathbb{F}_2 \}$ be a linear subspace of $\mathbb{F}_2^4$ and $\pi(x)=(x_1+\alpha x_2x_3+\beta x_2x_4+\gamma x_3x_4+\delta x_2x_3x_4,x_2,x_3 ,x_4):\alpha,\beta,\gamma,\delta \in \mathbb{F}_2 $ be a permutation of $\mathbb{F}_2^4$. Then $f_{(\pi, L)}:\mathbb{F}_2^4\times \mathbb{F}_2^4 \to \mathbb{F}_2$ is a $\mathcal{C}$-class bent function when $\gamma=\delta=0$ and $f_{(\pi, L)}$ is not a bent function otherwise.
\end{example}

Let us provide some justification for this. Clearly when $\gamma=\delta=0$, then by (B) of Theorem 4.1,  $(\pi, L)$ satisfies $\mathcal{C}$-condition. Let $D_{(\pi, L)}$ denote the support of $f_{(\pi, L)}$. When $\gamma =1$, for $\xi = (1,1,1,1,-1,1,1,1) \in U$, we see that $\big(\rho_{G}(D_{(\pi, L)})(\xi)\big)^2-\lambda \rho_{G}(G)(\xi)-(k-\lambda)=-64 \ne 0$. When $(\gamma, \delta) = (0, 1)$, $(\alpha, \beta) \ne (1, 1)$, for $\xi = (1,-1,1,1,-1,1,1,1)$ we have $\big(\rho_{G}(D_{(\pi, L)})(\xi)\big)^2-\lambda \rho_{G}(G)(\xi)-(k-\lambda)=-64 \ne 0$. Further when $\alpha, \beta, \gamma, \delta) = (1, 1, 0, 1)$, for $\xi = (1,-1,1,1,-1,1,1,1)$ we have $\big(\rho_{G}(D_{(\pi, L)})(\xi)\big)^2-\lambda \rho_{G}(G)(\xi)-(k-\lambda)=192 \ne 0$.

When $\alpha = \beta = \gamma =0$ and $\delta = 1$, we get a special case of the family in Section 5. Using Matlab, we have also determined that its Walsh-Hadamard spectrum (see \cite{MaG}) contains -16 with multiplicity 104, 16 with multiplicity 88, each of -32 and 32 with multiplicity 8 and 0 with multiplicity 48. It has also been verified that in the other non-bentness cases of Example 6.1, Walsh-Hadamard spectrum is not contained in $\{\pm 2^m\}$. This provides another verification of non-bentness, following the definition in \cite{MaG}.

\begin{example}
	Let $L=\{(0,0,x_3,x_4,x_5,x_6): x_3,x_4,x_5,x_6\in \mathbb{F}_2 \}$ be a linear subspace of $\mathbb{F}_2^6$ and $\pi(x)=(x_1+x_3x_4+x_5x_6,x_2,x_3 ,x_4,x_5,x_6)$ be a permutation of $\mathbb{F}_2^6$. Then $f_{(\pi, L)}:\mathbb{F}_2^6\times \mathbb{F}_2^6 \to \mathbb{F}_2$ is not a bent function.
\end{example}

This can be justified by observing that $\big(\rho_{G}(D_{(\pi, L)})(\xi)\big)^2-\lambda \rho_{G}(G)(\xi)-(k-\lambda)= -768 \ne 0$ where $D_{(\pi, L)}$ is the support of $f_{(\pi, L)}$ and $\xi = (1,1,1,1,1,1,-1,1,1,1,1,1) \in U$.
\\

\section{ Concluding Remarks}
 In this paper, we have connected the count of abelian difference sets with given parameters to computation of Hilbert functions of an ideal. There are algorithms for computation of Hilbert function. While diffculties in implementation for large values of parameters need to be addressed, the theoretical consequences of this connection can also be explored.

We also undertook to explore bentness of $f_{(\pi, L)}$ when $(\pi, L)$ does not satisfy $\mathcal{C}$-condition. Theorem 4.1 (C) helped us determine the choice of  $(\pi, L)$ for exploration and Sections 5 and 6 provided the results of exploration.This work was prompted by the following questions which still await the answers.

{\bf Question 1.} Is $\mathcal{C}$-condition necessary for bentness of $f_{(\pi, L)}$? If yes, then provide the proof or else provide the counter-example.

{\bf Question 2.} As a consequence of Theorem 4.1 (B) and (C), it follows that if $L=\{(\mathbf{0}_s,x_{s+1},\dots,x_{m}): x_i \in \mathbb{F}_2, s+1\le i \le m \}$ and $\pi_i(x_1, \dots, x_m)=(x_1+P_i(x_2, \dots, m_m), x_2, \dots, x_m)$ for $i=1, 2$ are such that $(\pi_i, L)$ satisfies $\mathcal{C}$-condition for $i=1, 2$ then $(\pi_1 \circ \pi_2, L)$ also satisfies $\mathcal{C}$-condition. What can we say about bentness of $f_{(\pi_1 \circ \pi_2, L)}$ if we know bentness of $f_{(\pi_i, L)}$ for $i=1,2$? In general, for a given subspace $L$ of $\mathbb{F}_2^m$, is there a semigroup structure on the set of all permutations $\pi$ of $\mathbb{F}_2^m$ such that $f_{(\pi, L)}$ is bent? If not, what are the counterexamples?

We hope to continue our exploration further guided by these questions.
\\

\noindent{\bf Acknowledgements :} Both the authors thank the support from FIST Programme vide SR/FST/MSI-090/2013 of DST, Govt. of India. The second author thanks UGC, Govt of India for the support under SRF Programme (SR. No. 2061540979, Ref. No. 21/06/2015(1)EU-V R. No. 426800).

\end{document}